\documentclass[12pt]{article}
\usepackage{amsmath}
\usepackage{amssymb}
\usepackage{amsthm}
\usepackage{amsfonts}
\usepackage{ulem}
\usepackage{cancel}
\usepackage{graphicx}
\usepackage{dsfont} 

\usepackage{bbm,epsfig,graphics,epic,color,rotating,color}
\numberwithin{equation}{section}

\newcommand{\R}{{\mathbb R}}
\newcommand{\Z}{{\mathbb Z}}

\newcommand{\N}{{\mathcal{N}}}
\newcommand{\Proof}{{\textit {Proof} }}

\newcommand{\ovl}[1]{\overline{#1}}
\newcommand{\wt}[1]{\widetilde{#1}}

\newcommand{\un}[1]{\underline{#1}}

\newcommand{\ed}{\mathrm{d}}

\newcommand{\reff}[1]{(\ref{#1})}

\newcommand{\E}{{\mathsf{E}}}
\renewcommand{\k}{\varkappa}
\newcommand\e{\varepsilon}
\newcommand\la{\lambda}
\renewcommand\phi{\varphi}

\newcommand{\be}{\begin{equation}}
\newcommand{\ee}{\end{equation}}
\newcommand{\bel}[1]{\begin{equation}\label{#1}}
\newcommand{\bea}{\begin{eqnarray}}
\newcommand{\eea}{\end{eqnarray}}
\newcommand{\balign}{\begin{align}}
\newcommand{\ealign}{\end{align}}
\newcommand{\ba}{\begin{array}}
\newcommand{\ea}{\end{array}}
\newcommand{\bfig}{\begin{figure}}
\newcommand{\efig}{\end{figure}}

\newcommand{\C}{{\mathbb C}}

\newtheorem{theorem}{Theorem}[section]
\newtheorem{lemma}[theorem]{Lemma}

\newtheorem{corollary}[theorem]{Corollary}

\newtheorem{remark}[theorem]{Remark}

\theoremstyle{definition}

\title{Large Fluctuations in two-level systems with stimulated emission}

\author{ E. Pechersky$^{1}$,\and S. Pirogov$^{1}$,\and G. M. Sch\"utz$^{2}$, \and A.~Vladimirov$^{1}$  and  A. Yambartsev$^3$ }
\date{}

\begin{document}

\maketitle

$^1$ Institute for Information Transmission Problems, 19, Bolshoj Karetny, Moscow, 127994, RF


$^2$Interdisziplin\"ares Zentrum f\"ur Komplexe Systeme, Universit\"at Bonn, Br\"uhler Str. 7, 53119 Bonn, Germany

$^3$Institute of Mathematics and Statistics, University of S\~ao Paulo (USP), S\~ao Paulo 05508-090,
SP, Brazil

\begin{abstract}
We consider a system of $N$ identical independent Markov processes, each taking values 0 or 1. The system describes a stochastic dynamics of an ensemble of two-level atoms. The atoms are exposed to a photon flux.  Under the photon flux action, every atom changes its state with some intensities either from its ground state (state 0) to the excited state (state 1) or from the excited state to the ground state (the stimulated emission). The atom can also change the state from the excited  to the ground one spontaneously. 

We study rare events when the big cumulative emission occurs on the fixed time interval $[0,T]$. To this end we apply the large deviation theory which allows one to make asymptotic analysis as $N\to\infty$.
\end{abstract}



\section{Introduction} This article studies a generalisation of the model introduced in \cite{PP}. Some of the results obtained in this article are similar to those in  \cite{PP}. The model considered here contains a new feature which makes the model closer to  known models resembling the studies of A. Einstein on the interactions of the matter and the light (see, for example, \cite{WH}).

Informally the model in \cite{PP} is constituted by a set of ``atoms" and a flux of ``photons"  falling on the atoms. It is assumed that the atoms have two levels of energy $E_1$  and $E_2$, $E_1<E_2$, two-level atoms. The level $E_1$ is called \textit{ground state} and the level $E_2$ is called \textit{excited state}. The stochastic evolution in this model is formed by jumps of the atom states between the energy levels. The jumps of the atoms from the ground states to the excited one happen with  rate $\lambda$  under the action of the photon flux. The inverse jumps happen spontaneously with  rate $\mu$. There are no restrictions either for  values  $\lambda$, $\mu$  or for relations between these parameters. 

The goal of \cite{PP} was to study the large  cumulative emission of the photons during the time interval $[0,T]$. Such type of the problems can be solved by a part of the probability theory called the large deviation theory. The main tool of the theory for solving such problems is a \textit{rate function}. The large deviation theory gives asymptotic answers about the probabilities of the large fluctuations under some proper scaling of the studied random system. Such analysis was performed in \cite{PP}. 

In this article we study a similar system. On contrast to the system in \cite{PP}, we add a new process here:   a \textit{stimulated emission}. The stimulated emission occurs when a photon from the flux knocks out another photon from an atom  in the excited state.  Then the atom emits two photons. Our interest is in the large emission during the time interval $[0,T]$,  it is the same as in \cite{PP}.  The problem here is  also a subject of the large deviation theory. The rate function of the studied model is different from   that in \cite{PP}; however, the forms of both  rate functions are similar.

It is worth to mention that the rate function contains a lot of information about properties of large fluctuations. In particular, it includes  information on the average path  to achieve these fluctuations.  The extraction of this information is reduced to a solution of a Hamiltonian system of  differential equations. This work was done in \cite{PP} where the Hamiltonian system consists of four equations with four unknown functions on $[0,T]$, and the two of the equations depend on two functions only. This makes it possible to to solve full system explicitly.

In the case considered here, the situation is quite similar. There are six differential equations, three of which depend only on three unknown functions. This allows one, as in \cite{PP}, to find solutions. The solutions describe properties of the path leading to the large fluctuations.
The main result of the work is a behaviour of the Markov process which attains large emission. It was established that given the large emission on $[0,T]$ the ratio of the numbers of the excited atoms and the atoms in the ground state is close to 1, the greater the emission the closer the ratio to 1 (see Section~\ref{cha}).
This also holds for the model in \cite{PP}.

However, the presence of the stimulated emission causes a new effect which is not present in the model from \cite{PP}. The new property consists in the fact that for given large fluctuation of the emission the balance between the number of photons emitted by the stimulations and the photons emitted spontaneously is shifted to the stimulated photons. In the limit when the emission tends to infinity only stimulated photons give
 contribution to the radiation. 

The effect looks rather unexpected. The corresponding theorem is proved in Section \ref{las}.

\section{Model and generators} A formal representation of the model is the Markov process $\Xi(t)=(\Xi_1(t),\Xi_2(t),\Xi_3(t))$ taking values in three-dimensional space  $\Z_+^3$.

Firstly, we define $\Xi_1(t)$ which itself is the Markov process. To this end, consider a space configuration $\cal X$ defined on a finite set $\N=\{1,2,...,N\}$ and taking their values in $\{0,1\}$,  $\mathcal{X}= \{0,1\}^{\N}$. The set $\N$ represents an ensemble of two-level atoms and a configuration $\un x\in\cal X$ is a map
\[\un x:\:\N\to \{0,1\}\]
describing the states of all atoms in $\N$. 
The value $\un x(i)=0$ means that the atom $i\in \N$ is in the ground state, and $\un x(i)=1$ means that the atom is excited. $\Xi_1$ is the number of excited atoms in the system. Thus, the values which the process $\Xi_1$ takes are
\[M=M_{\un x}=\sum_{i\in\N}\un x(i).\]
Atoms are excited with rate $\lambda$, corresponding to a jump $M \to M+1$, and return to
the ground state with rate $\tilde{\mu}$, corresponding to a jump $M \to M-1$. 
Transitions of 
$\Xi_1$ are therefore determined by the infinitesimal generator.
\begin{equation}\label{2.1}
\mathbf{\tilde{L}}G(M)=
\lambda(N- M)\left[G(M+1)- G(M) \right]+\wt\mu M\left[ G(M-1)- G(M)\right],
\end{equation}
The domain of $\mathbf{\tilde{L}}$ is the set of all real functions $G$ on $\Z_+$. 

The process $\Xi_1$ has jump-wise  paths. The values of the jumps are $+1$ or $-1$. If $\Xi_1(t)=M\ge 0,$ then the expected time till the next jump is $(\la(N-M)+\wt\mu M)^{-1}$.  Let $t'$ be the time when this jump happens then 
\begin{enumerate}
\item[]
$\Xi_1(t'+)-\Xi_1(t)=+1$ with the probability $\frac{(N-M)\la}{(N-M)\la+\wt\mu M}$ and  
\item[] $\Xi_1(t'+)-\Xi_1(t)=-1$ with the probability $\frac{\wt\mu M}{(N-M)\la+\wt\mu M}$.
\end{enumerate}

The two components $\Xi_2(t)$ and $\Xi_3(t)$ of the process $\Xi(t)$ describe the number of emitted photons on  time interval $[0,t]$. It is always assumed that $\Xi_2(0)=\Xi_3(0)=0$. The process $\Xi_2(t)$ is the number of spontaneously emitted photons while $\Xi_3(t)$ is the number of the stimulated emitted photons. 
\begin{enumerate}
\item[]
If $\Xi_1(t)>0$ at the time $t$ and the spontaneous emission occurs then $\Xi_2(t+)-\Xi_2(t)=1$, and at this moment $\Xi_1(t+)-\Xi_1(t)=-1$, and $\Xi_3(t)$ does not change its value. 
\item[]
If $\Xi_1(t)>0$ at the time $t$ and the stimulated emission occurs then $\Xi_3(t+)-\Xi_3(t)=2$, two photons are emitted. At this moment $\Xi_1(t+)-\Xi_1(t)=-1$ and $\Xi_2(t)$ does not change its value. 
\end{enumerate}
The intensity of spontaneous emission is $\mu$ and the intensity of stimulated emission is $\nu$,
$\mu+\nu=\tilde{\mu}$ (see \reff{2.1}).
The stochastic dynamics of $\Xi(t)=(\Xi_1(t),\Xi_2(t),\Xi_3(t))$ is driven by the infinitesimal operator 
\[
\begin{array}{ll}
\mathbf{L}G(M_1,M_2,M_3)=&\lambda(N- M_1)\left[G(M_1+1,M_2,M_3)- G(M_1,M_2,M_3) \right]\\
&{} + \mu M_1\left[ G(M_1-1,M_2+1,M_3)- G(M_1,M_2,M_3)\right]\\
&{} + \nu M_1\left[ G(M_1-1,M_2,M_3+2)- G(M_1,M_2,M_3)\right].
\end{array}
\]
The domain of the functions $G$ is $\Z_+^3$.

\section{Scaling} In this section we consider the scaling of  $\Xi(t)$ by the number of atoms $N$. We study the asymptotic of the system as $N\to\infty$. The  scaled process is
\[
\xi^N(t)\equiv (\xi^N_1(t),\xi^N_2(t),\xi^N_3(t)) :=\frac1N\Xi(t)\equiv \left(\frac1N\Xi_1(t),\frac1N\Xi_2(t),\frac1N\Xi_3(t)\right)
\]
The process $\xi^N$ takes its value in $\frac1N \Z_+^3$. The set of jumps is 
\[
\left\{\left(+\frac1N,0,0\right),\left(-\frac1N,+\frac1N,0\right)\left(-\frac1N,0,+\frac2N\right)
\right\}
\]
The generator of $\xi^N$ is
\[
\begin{array}{ll}
\mathbf{L}_Ng(x_1,x_2,x_3)=&N\lambda(1-x_1)\left[g\left(x_1+\frac1N,x_2,x_3\right)- g(x_1,x_2,x_3) \right]\\
&{} + N\mu x_1\left[ g\left(x_1-\frac1N,x_2+\frac1N,x_3\right)- g(x_1,x_2,x_3)\right]\\
&{} + N\nu x_1\left[ g\left(x_1-\frac1N,x_2,x_3+\frac2N\right)- g(x_1,x_2,x_3)\right],
\end{array}
\]
where $(x_1,x_2,x_3)\in\frac1N\Z_+^3\subset\R_+^3$.

Let $\mathbb F_N$ be the set of all bounded function on $\frac1N\Z_+^3$ and let $P_N:\:\C_0^1(\R_+^3)\to \mathbb F_N$ is the restriction operator.  

Define
\[
\begin{array}{ll}
\mathbf{L}_\infty g(x_1,x_2,x_3)=&\lambda(1-x_1)\frac{\partial g(x_1,x_2,x_3)}{\partial x_1}+\mu x_1\left(\frac{\partial g(x_1,x_2,x_3)}{\partial x_2}-\frac{\partial g(x_1,x_2,x_3)}{\partial x_1}\right)\\
& {} + \nu x_1\left(2\frac{\partial g(x_1,x_2,x_3)}{\partial x_3}-\frac{\partial g(x_1,x_2,x_3)}{\partial x_1}\right).
\end{array}
\]
It is clear that
\begin{equation}\label{aa}
\lim_{N\to\infty}\sup_{(x_1,x_2,x_3)\in\frac1N\Z_+^3}\bigl|\mathbf{L}_\infty P_Ng(x_1,x_2,x_3)-\mathbf{L}_Ng(x_1,x_2,x_3)\bigr|=0
\end{equation}
for any $g\in\C_0^1(\R_+^3)$.  

Having (\ref{aa}) the following result is a consequence of the Trotter-Kurtz
theorem (\cite{Ku} or \cite{EK},Ch.1).

\begin{theorem}\label{t3} Let $\xi_N(0)$ tend to some point $(u_1^0,u_2^0,u_3^0)\in\R_+^3$, then for any $t>0$ the random vector $\xi_N(t)$ tends in probability to the non-random limit $(x_1(t),x_2(t),x_3(t))$. The functions  $(x_1(t),x_2(t),x_3(t))$ are the solution of the following system
\begin{equation}\label{3.01}
\begin{array}{lll}
\dot x_1(t)&=&\la(1-x_1(t))-(\mu+\nu) x_1(t)\\
\dot x_2(t)&=&\mu x_1(t)\\
\dot x_3(t)&=&\nu x_1(t)
\end{array}
\end{equation}
with the initial conditions $x_1(0)=u_1^0,x_2(0)=u_2^0,x_3(0)=u_0^3$.
\end{theorem}

\section{Problems} The first problem we address in this work is to find an ``optimal" path of the emission when a large emission occurs. The emission path on $[0,T]$ is the trajectory of $\Xi_2(t)+\Xi_3(t),t\in [0,T]$. As we mentioned above $\Xi_2(0)+\Xi_3(0)=0$. Let  ${\cal A}(A)=(\Xi_2(T)+\Xi_3(T)\geq A)$ be an event meaning that  the total emission at the time $T$ is at least $A$. The question is: what is the asymptotic of the logarithm of the probability $\ln\Pr({\cal A}(A))$ as $A\to\infty$? This asymptotic can be  equivalently expressed in terms of the processes $\xi_N$. Let $A=BN$ then $\Pr({\cal A}(A))=\Pr(\xi_2^N(T)+\xi_3^N(T)\geq B)$. 
The problem of optimal path is to find $\E\left((\xi_2^N(t)+\xi_3^N(t))/(\xi_2^N(T)+\xi_3^N(T))\geq B\right)$, which is a function of $t\in[0,T]$, and its limit as $N\to\infty$.

The second problem is to find behaviour of the function $\E\left(\xi_1^N(t)/(\xi_2^N(T)+\xi_3^N(T))\geq B\right)$ for large $N$.

Relations between $\xi_2^N(T)$ and $\xi_3^N(T)$ on the event $(\xi_2^N(T)+\xi_3^N(T)\geq B)$ is the third problem.

\section{Large deviations} The main tool we use is the large deviation theory. We can apply the large deviation theory since Theorem \ref{t3}. The theorem establishes the convergence of the process $\xi^N$ to a non-random path. The pre-limiting processes $\xi^N$ belong to the set $D[0,T]$ of step-wise trajectories. Therefore we introduce Skorohod metric in $D[0,T]$ to enable the large deviation principle.
The main tool of the large deviation theory is the \textit{rate function} (see \cite{DZ}, \cite{FeK}.)

For functions $x_1, x_2, x_3 \in\C^1[0,T]$ (differentiable on interval $[0,T]$) with condition $x_2(0)=x_3(0)=0$, the rate function in our case is
\begin{eqnarray}\label{rate}
I(x_1,x_2,x_3)&=&\int_0^T \sup_{\varkappa_1(t),\varkappa_2(t),\varkappa_3(t)}\bigl(\varkappa_1(t)\dot x_1(t)+\varkappa_2(t)\dot x_2(t)+\varkappa_3(t)\dot x_3(t)-\\
&& \ \ - \rho(x_1(t))[\varphi(\varkappa_1(t),\varkappa_2(t),\varkappa_3(t))-1]\bigr)\ed t, \nonumber
\end{eqnarray}
where $\rho(x_1(t))=\lambda(1- x_1(t))+\mu x_1(t)+\nu x_1(t)$, and
\begin{equation}\label{laplace}
\varphi(\varkappa_1,\varkappa_2,\varkappa_3)=\frac{\lambda (1-x_1)}{\rho(x_1)} e^{\varkappa_1}+\frac{\mu x_1}{\rho(x_1)} e^{- \varkappa_1 + \varkappa_2}+\frac{\nu x_1}{\rho(x_1)} e^{- \varkappa_1 + 2\varkappa_3},
\end{equation}
i.e.
\begin{eqnarray}
I(x_1,x_2,x_3) &=&
 \int_0^T \sup_{\varkappa_1,\varkappa_2,\varkappa_3}\bigl(\varkappa_1\dot x_1+\varkappa_2\dot x_2+\varkappa_3\dot x_3-\lambda(1-x_1)[e^{\varkappa_1}-1] \label{4.1}\\
 &&{} - \mu x_1[e^{-\varkappa_1+\varkappa_2}-1]-\nu x_1[e^{-\varkappa_1+2\varkappa_3}-1] \bigr)\ed t.\nonumber
\end{eqnarray}

The rate function in our case is defined on the space $\C_{[0,T]}(\R_+^3)$. The value of the rate functions for the functions having discontinuity is equal to infinity.

The expression \reff{4.1} for the rate function is derived from the Sanov theorem and the contraction principle (see \cite{DZ}, \cite{Reza}). Also for short explanations see Remark~2.1 in \cite{PP}. Note that we can apply the Sanov theorem because the stochastic dynamics of atoms are independent from each  other.

Let $Z_N=\{(x_1,x_2,x_3):\:x_2(0)=x_3(0)=0,x_2(T)+x_3(T)\geq B\}\subset(\C^1(0,T])^3$. In order to find the optimal path we have to find a path $(x_1^0,x^0_2,x^0_3)$ from $Z_N$ where $\inf_{(x_1,x_2,x_3)\in Z_N} I(x_1,x_2,x_3)$ is attained: 
$$
I(x_1^0,x^0_2,x^0_3)=\inf_{(x_1,x_2,x_3)\in Z_N} I(x_1,x_2,x_3).  
$$
To this end we have to solve the following system of equations
\begin{equation}\label{hameq}
\left\{ \begin{array}{rcl}
\dot x_1&=&\lambda (1-x_1)\exp\{\varkappa_1\}+\mu x_1\exp\{-(\varkappa_1-\varkappa_2)\}+\nu x_1\exp\{-(\varkappa_1-2\varkappa_3)\}, \\
\dot x_2&=&\mu x_1\exp\{-(\varkappa_1-\varkappa_2)\}, \\
\dot x_3&=&2\nu x_1\exp\{-(\varkappa_1-2\varkappa_3)\}, \\
\dot \varkappa_1&=&\lambda \exp\{\varkappa_1\}-\mu\exp\{-(\varkappa_1-\varkappa_2)\} -\nu\exp\{-(\varkappa_1-2\varkappa_3)\}- \lambda + \mu+\nu, \\
\dot\varkappa_2&=&0\\
\dot\varkappa_3&=&0
\end{array} \right.
\end{equation}
which can be considered as a Hamiltonian system generated by the Hamiltonian
\[
H(\ovl x,\ovl \k)=\lambda(1-x_1) \bigl[e^{\varkappa_1}-1 \bigr]-
 \mu x_1 \bigl[e^{-\varkappa_1+\varkappa_2}-1 \bigr]-\nu x_1 \bigl[e^{-\varkappa_1+2\varkappa_3}-1\bigr] ,
\]
where $\ovl x=(x_1,x_2,x_3), \ovl \k=(\k_1,\k_2,\k_3)$ and $x_i, \k_i \in\C^1[0,T], i=1, 2, 3$.

A peculiarity of the system is that the last three equations  do not depend on the paths $x_1,x_2,x_3$. It allows one to solve the fourth equation.  
The method of the solution is the same as in \cite{PP}. We get 
\begin{equation}\label{sol2}
e^{\varkappa_1 (t) } = \frac{r_2-r_1C_1\exp\{t\lambda(r_2-r_1)\}}{1-C_1\exp\{t\lambda(r_2-r_1)\}},
\end{equation}
where
\begin{equation}\label{r12}
r_{1,2} = \frac{a}{2} \mp \sqrt{ \left( \frac{a}{2} \right)^2 + b },
\end{equation}
\begin{equation}
\label{notation}
\begin{array}{l}
\gamma=\frac{\mu}{\lambda},\\[0.3cm]
\gamma_1=\frac{\nu}{\lambda},\\[0.3cm]
a=1-\gamma-\gamma_1,\\[0.3cm]
b=\gamma e^{\varkappa_2}+\gamma_1 e^{2\varkappa_3},
\end{array}
\end{equation}
$ C_1 $ is a constant which can be found from the boundary condition  $\varkappa_1(T)=0$. The last equality holds  since there are no constraints on the value $x_1(T)$. Thus 
\begin{equation}\label{C1}
C_1=\frac{r_2-1}{r_1-1}\exp\{-T\lambda(r_2-r_1)\}.
\end{equation}
We can find the path 
 \begin{equation}\label{f1}
x_1(t)=\frac{\lambda(e^{\varkappa_1}-1)-K}{\lambda(e^{\varkappa_1}-1)+\mu(1-e^{-\varkappa_1+\varkappa_2})+\nu(1-e^{-\varkappa_1+2\varkappa_3})},
\end{equation}
where $K$ is a constant which can be found from the above equality at $t=0$ considering $x_1(0)$ as a parameter:
\begin{equation}\label{K}
K=\lambda(e^{\varkappa_1(0)}-1)(1-x_1(0))-x_1(0)\mu(1-e^{-\varkappa_1(0)+\varkappa_2})+x_1(0)\nu(1-e^{-\varkappa_1(0)+2\varkappa_3}).
\end{equation}
Now we find the functions
\begin{equation}\label{x23}
\begin{aligned}
x_2(t) &= \mu e^{\varkappa_2}\int\limits_0^tx_1(s)e^{-\varkappa_1(s)}\ed s,\ \ t\in[0,T], \\
x_3(t) &= 2\nu e^{2\varkappa_3}\int_0^tx_1(s)e^{-\varkappa_1(s)}\ed s,\ \ t\in[0,T].
\end{aligned}
\end{equation}
\section{Emerging chaos}\label{cha} In this section we study  the behaviour of the Markov process with  large $B=x_2(T)+x_3(T)$.  The following theorem shows that the system approaches to the maximal chaos as $B\to\infty$.
\begin{theorem}\label{t1}
For any $0<\alpha<T/2$ and any $\e>0$, there exists $B_{0}\equiv B_{0}(\alpha, \e)>0$ such that
the inequality
\[
\Bigl| x_{1}(t)-\frac{1}{2} \Bigr| <\e
\]
holds for all $t\in[\alpha,T-\alpha]$ whenever $B\ge B_{0}$.
\end{theorem}
\begin{remark}
The theorem claims that at the large emission limit, as $B\to\infty$, in  the system of  two-level atoms the probability of the atoms  to have  their states 0 or 1 tends to $\frac12$. 
\end{remark}

The proof basically repeats the proof of the similar result in \cite{PP} with some modifications. 
\begin{lemma}
The following relation
\begin{equation}\label{eess}
r_2\geq cB  
\end{equation}
holds for some $0<c<1$ if $B$ is large enough.
\end{lemma}
\proof 
We consider two cases: $r_2>1$ and $r_2\leq 1$. 

In the case  $r_2>1$ 
\begin{equation}\label{kaper}
e^{\varkappa_1(t)}\leq r_2,
\end{equation}
using directly \reff{sol2} and \reff{C1} and noting that always $r_1 < 0 < r_2$. Further, because of $\varkappa_1(T)=0$, by integrating fourth equation of \reff{hameq} we obtain
\begin{equation}\label{kap1}
\varkappa_{1}(0)=\left(\mu e^{\k_2}+\nu e^{2\k_3}\right)\int_{0}^{T}e^{-\k_1(t)}dt-\lambda\int_{0}^{T}e^{\k_1(t)}dt+(\lambda-\mu-\nu)T.
\end{equation}
Condition $x_2(T)+x_3(T)\ge B$ using \reff{x23} gives us 
\begin{equation}\label{kap2}
\left(\mu 
e^{\varkappa_2}+2\nu e^{2\varkappa_3}\right)\int_0^Tx_1(t)e^{-\varkappa_1(t)}\ed t \ge B.
\end{equation}
The inequality
\[
\frac{\left(\mu 
e^{\varkappa_2}+\nu e^{2\varkappa_3}\right)\int_0^Te^{-\varkappa_1(t)}\ed t}{\int_0^Tx_1(t)e^{-\varkappa_1(t)}\ed t}\geq \mu 
e^{\varkappa_2}+\nu e^{2\varkappa_3}
\]
is obvious since $0\leq x_1\leq 1$. Then
\begin{equation*}
\k_1(0)\geq \left(\mu 
e^{\varkappa_2}+\nu e^{2\varkappa_3}\right)\int_0^Tx_1(t)e^{-\varkappa_1(t)}\ed t -\lambda\int_{0}^{T}e^{\k_1(t)}dt+(\lambda-\mu-\nu)T,
\end{equation*}
which implies with \reff{kap2}
\[
\k_1(0)+\lambda\int_{0}^{T}e^{\k_1(t)}dt\geq B+(\lambda-\mu-\nu)T.
\]
Using now inequality \reff{kaper} we obtain
\[
\ln r_2 + \la Tr_2 \geq B+(\lambda-\mu-\nu)T.
\]
Let $s= \la Tr_2$ and $S=B+(\lambda-\mu-\nu)T+\ln(\la T)$. Then the above inequality is
\begin{equation}\label{eq111}
\ln s + s \geq S.
\end{equation}
We prove that inequality \reff{eq111} implies 
\[
s\geq S-\ln S
\]
for $S>1$.
Indeed, we can assume that $S>1$ since $B$ is large. If $s>S$ then $\ln s>\ln S>0$ and 
\[
s\geq S\geq S-\ln S.
\]
If $s\leq S$ then, from \reff{eq111}
\[
s\geq S-\ln s\geq S-\ln S.
\]
Thus, for any  $0<c'<1$,
\[
\la Tr_2>c'B
\]
for $B>\frac1{1-c'}(\ln(B+R)-R)$, where $R=(\la-\mu-\nu)T+\ln(\la T)$

The second case is $r_2<1$ for large $B$. In this case $\sup_Bb<\infty$ and
\[
\sup_B\sup_{t\in[0,T]}e^{-\k_1(t)}<\infty
\]
which contradict to \reff{kap2}. \qed

\begin{corollary}
\begin{equation}\label{7.0}
\lim_{B\to\infty}(\mu e^{\k_2}+\nu e^{2\k_3})\to\infty
\end{equation}
as $B\to\infty$
\end{corollary}
The proof follows from \reff{r12} and \reff{notation} \qed
\begin{corollary}
\begin{equation}\label{7.1}
\begin{aligned}
e^{\k_1}=r_2+o\bigl( (\mu e^{\varkappa_2}+\nu e^{2\varkappa_3})^{\frac12} \bigr)&=\Bigl( \frac\mu\lambda e^{\varkappa_2}+\frac\nu\lambda e^{2\varkappa_3} \Bigr)^{\frac12}+o\bigl( (\mu e^{\varkappa_2}+\nu e^{2\varkappa_3})^{\frac12} \bigr)\\
e^{-\k_1} \Bigl( \frac\mu\lambda e^{\varkappa_2} + \frac\nu\lambda e^{2\varkappa_3} \Bigr)&=\Bigl( \frac\mu\lambda e^{\varkappa_2}+\frac\nu\lambda e^{2\varkappa_3} \Bigr)^{\frac12}+o\bigl( (\mu e^{\varkappa_2}+\nu e^{2\varkappa_3})^{\frac12} \bigr)\
\end{aligned}
\end{equation}
uniformly for $t\in[0,T-\alpha)$, where $\alpha<\frac T2$.

\end{corollary}
\proof For any $t\in[0,T-\alpha)$ 
\[
\lim_{B\to\infty}|r_1C_1|\exp\{t\la(r_2-r_1)\}=0.
\]
Then the relations \reff{7.1} follow from \reff{sol2}. \qed

The corollary proves the following lemma.

\begin{lemma}\label{L1}
The limit
\begin{equation}\label{asym}
\frac{\la\exp\{2\k_1(t)\}}{\mu e^{\varkappa_2}+\nu e^{2\varkappa_3}}\to 1
\end{equation}
as $B\to\infty$, holds true uniformly on $t\in[\alpha,T-\alpha]$ for any fixed $\alpha<\frac T2$ . 
\end{lemma} 
\vspace{0.5cm}

\Proof of Theorem~\ref{t1}. The proof repeats the similar consideration in \cite{PP} with some modifications.
Having \reff{asym}  we rewrite the equation for $\dot x_1$ 
(see \reff{hameq}) on the interval $[0,T-\alpha]$ as follows
\be\label{f1ap}
\dot x_1(t)=\la\exp\{\k_1(t)\}\big(1-x_1(t)\big)-\la\exp\{\k_1(t)\}(1+\delta(t))x_1(t),
\ee
where $\delta(t)\to 0$  uniformly on $[0,T-\alpha]$ as $B\to\infty$.  Equivalently, the equation above is
\begin{equation}\label{f1ap2}
\dot x_1(t)=\la r_2(1-2x_1(t)) +\varepsilon_0(t),
\end{equation}
where $|\e_{0}(t)|\le\e$ for all $t\in[0,T-\alpha]$. Let $g(t)=x_1(t)-\frac12$. Then 
\begin{equation}\label{f1ap3}
\dot g(t)=\la r_2(1+2g(t)) +\varepsilon_0(t).
\end{equation}
The solution of the equation on $[0,T-\alpha]$ is
\begin{equation}
g(t)=e^{-2\la r_2}g(0)+\int_0^{\la r_2t}\e_0(s)e^{2\la r_2(s-t)}\ed s\to 0
\end{equation}
as $B\to\infty$. Hence,  
$$
x_1(t)\to\frac12
$$
on $[0,T-\alpha]$.
\qed

\section{Balance between the spontaneous and stimulated emissions.}\label{las} This section is devoted to the studies of the relation between fractions of spontaneous and stimulated emissions subject to the total large emission. The next theorem claims a rather unexpected result in this respect.

Recall that the total emission is
\begin{equation}\label{5.3}
B=x_2(T)+x_3(T)=(\mu e^{\k_2}+2\nu e^{2\k_3})\int_0^Tx_1(s) e^{-\k_1(s)}\ed s
\end{equation}
(see \reff{x23}).
We used the fact that $\k_2$ and $\k_3$ are constant over  time. 
\begin{theorem}\label{bal}
\[
\lim_{B\to\infty}\frac{x_2(T)}{x_3(T)}=0.
\]
\end{theorem}
\proof We have to minimise the rate function $I$ \reff{4.1} over $x_2$ and $x_3$ given large $B$, (see \reff{5.3}). To this end we find the limit of every term in \reff{4.1} as $B\to\infty$.
Using \reff{r12} we obtain
\begin{equation}\label{8.0}
r_{1,2}=\mp\sqrt{b}\phi_{\mp} \cong \mp\sqrt{b},
\end{equation}
as $b\to\infty$, which is the same as $B\to\infty$, where  
\begin{equation}\label{9.0}
\phi_{\mp}:=\mp\sqrt{1+\frac1b\left(\frac a2\right)^2 }+\frac1{\sqrt{b}}\left(\frac a2\right).
\end{equation} 
Let 
\begin{equation}\label{9.00}
\begin{array}{ll}
\psi:=\sqrt{1+\frac1b\left(\frac a2\right)^2 },\\[.3cm]
\zeta_1:=\frac{r_2-1}{1-r_1}=\frac{\sqrt{1+\frac1b(\frac{a}{2})^2}+\frac1{\sqrt{b}}(\frac a2-1)}{\sqrt{1+\frac1b(\frac{a}{2})^2}-\frac1{\sqrt{b}}(\frac a2-1)},\\[.5cm]
\zeta_2:=-\frac{r_1}{r_2}=\frac{\sqrt{1+\frac1b(\frac{a}{2})^2}-\frac a{2\sqrt{b}}}{\sqrt{1+\frac1b(\frac{a}{2})^2}+
\frac a{2\sqrt{b}}}.
\end{array} 
\end{equation}
Using the expression:
\[
\begin{aligned}
-C_1\exp\{\la t(r_2-r_1)\}&=\frac{r_2-1}{1-r_1}\exp\{\la (t-T)(r_2-r_1)\}=\zeta_1\exp\left\{2\la\sqrt{b}\psi(t-T)\right\},\\
-\frac{r_1}{r_2}C_1\exp\{\la t(r_2-r_1)\} &=-\frac{r_1}{r_2}\frac{r_2-1}{1-r_1}\exp\{\la (t-T)(r_2-r_1)\}=
\zeta_1\zeta_2\exp\left\{2\la\sqrt{b} \psi(t-T)\right\}.
\end{aligned}
\]
we obtain 
\begin{equation}\label{8.1}
\begin{aligned}
e^{-\k_1}&=\frac{1-C_1\exp\{\la t(r_2-r_1)\}}{r_2-r_1C_1\exp\{\la t(r_2-r_1)\}}=
\frac1{r_2} \cdot \frac{1-C_1\exp\{\la t(r_2-r_1)\}}{1-\frac{r_1}{r_2}C_1\exp\{\la t(r_2-r_1)\}} \\ &=
\frac1{\sqrt{b}\phi_+} \cdot \frac{1+\zeta_1\exp\left\{2\la\sqrt{b}\psi(t-T)\right\}}{1+\zeta_1\zeta_2\exp\left\{2\la\sqrt{b} \psi(t-T)\right\}}.
\end{aligned}
\end{equation}
Some obvious limits as $b\to \infty$ are
\begin{equation}\label{limi}
\begin{aligned}
&\phi_-\to -1,\phi_+\to 1,\\
&\psi\to 1,\\
&\zeta_1\to 1,\zeta_2\to 1.
\end{aligned}
\end{equation}
Let 
\begin{equation}\label{8.3}
\begin{aligned}
\mu e^{\k_2}&=\alpha B^2,\\
2\nu e^{2\k_3}&=\beta B^2,
\end{aligned}
\end{equation}
where the relation between $\alpha$ and $\beta$ can be found from the equality
\begin{equation}\label{9.1}
B=(\mu e^{\k_2}+2\nu e^{2\k_3})\int_0^Tx_1(s)e^{-\k_1(s)}\ed s=B^2(\alpha+\beta)\int_0^Tx_1(s)e^{-\k_1(s)}\ed s.
\end{equation}
Evaluation of $\int_0^Tx_1(s)e^{-\k_1(s)}\ed s$ gives 
\begin{equation}\label{9.2}
b=\frac1\la(\mu e^{\k_2}+\nu e^{2\k_3})=\frac1\la B^2\left(\alpha+\frac{\beta}{2}\right).
\end{equation}

In what follows the arrow $\to$ means the limit as $b\to\infty$.

It follows from \reff{8.1} that
\[
\sqrt{b}\int_0^Tx_1(s)e^{-\k_1(s)}\ed s\to \int_0^Tx_1(s)\ed s=T\wt{x}_1,
\]
where $\wt{x}_1=\frac1T\int_0^Tx_1(s)\ed s$.

We derive the relation for $B$ from \reff{9.2} and plug it in into \reff{9.1}. We get
\begin{equation}\label{10.2}
\frac{\sqrt{\alpha+\frac\beta 2}}{\alpha+\beta}=\sqrt{\la}T\wt{x}_1(1 + F(b)),
\end{equation}
where $F(b)\to 0$, as $b\to\infty$.

Using the asymptotics which were found above we find the asymptotics of the rate function $I$ \reff{4.1} as $B\to \infty$. Assume that the functions $x_1,x_2,x_3,\k_1,\k_2,\k_3$ satisfy the equation system \reff{hameq}. We evaluate the asymptotics of $I$ separately calculating the terms in \reff{4.1}. Consider the first term $\int_0^T\k_1(s)\dot x_1(s) \ed s$. We use \reff{8.1} and the first equation in \reff{hameq}. So,
\begin{equation}\label{10.3}
\begin{aligned}
&\int_0^T\k_1(s)\dot x_1(s) \ed s=\\
&{} = \la\int_0^T\k_1(s)e^{\k_1(s)}\ed s-\la\int_0^Tx_1(s)\k_1(s)e^{\k_1(s)}\ed s+(\mu e^{\k_2}+\nu e^{2\k_3})\int_0^Tx_1(s)\k_1(s)e^{-\k_1(s)}\ed s.
\end{aligned}
\end{equation}
For the first two integrals in the right-hand side of \reff{10.3}, taking in account \reff{8.1} with limits \reff{limi}, and \reff{9.2}, we use the following asymptotics: 
$$
e^{\k_1} \cong \sqrt{b} = \frac{B}{\sqrt{\la}} \sqrt{\alpha+\frac{\beta}{2}} \mbox{ and } \k_1 \cong \log B.
$$
Thus,
\[
\frac{\la\int_0^T\k_1(s)e^{\k_1(s)}\ed s}{B\ln{B}}\to \sqrt{\la} \sqrt{\alpha+\frac{\beta}{2} }\cdot T
\]
for the first term in \reff{10.3},  and, using the same argument, we get 
\[
\frac{\la\int_0^Tx_1(s)\k_1(s)e^{\k_1(s)}\ed s}{B\ln{B}}\to\sqrt{\la}\sqrt{\alpha+\frac{\beta}{2} }\cdot T\wt{x}_1
\]
for the second one. 
The third term in \reff{10.3} has the limit
\[
\frac{(\mu e^{\k_2}+\nu e^{2\k_3})\int_0^Tx_1(s)\k_1(s)e^{-\k_1(s)}\ed s}{B\ln{B}}\to\sqrt{\la}\sqrt{\alpha+\frac{\beta}{2} }\cdot T\wt{x}_1.
\]

Next we find the asymptotics of the following two terms in $I$ (see \reff{4.1}). It is 
\[
\begin{aligned}
& \frac1{B\ln{B}}\left(\int_0^T\k_2\dot x_2(s) \ed s+\int_0^T\k_3\dot x_3(s) \ed s\right)=\\ 
&{} = \frac1{B\ln{B}}\left((\k_2\mu e^{\k_2}+2\k_3\nu e^{2\k_3})
\int_0^T x_1(s)e^{-\k_1(s)} \ed s\right)\to 2 \sqrt{\la}T\wt{x}_1 \sqrt{\alpha+\frac{\beta}{2}}.
\end{aligned}
\]
 
Now we note that the last terms of $I$
$$
- \int_0^T \Bigl( \lambda (1-x_1) \bigl( e^{\k_1} -1 \bigr) + \mu x_1 \bigl( e^{-\k_1+\k_2} -1 \bigr) +
\nu x_1 \bigl( e^{-\k_1+2\k_3} -1 \bigr) \Bigr) ds 
$$ 
divided by $B\ln{B}$, tends to 0 as $B\to\infty$.

Using all these calculations we have 
\[
\lim_{B\to\infty}\frac 1{B\ln{B}}I\to J(\alpha,\beta)=  (\sqrt{\la}T+4\sqrt{\la}T\wt{x}_1)\sqrt{\alpha+\frac{\beta}{2}}.
\]
The expression in the brackets grows with growing $\alpha$ for any $\beta$. We minimise $J(\alpha, \beta)$ over positive parameters $\alpha$ and $\beta$ under the condition \reff{10.2}, which in the limit takes the form 
\begin{equation}\label{10.22}
\frac{\sqrt{\alpha+\frac\beta 2}}{\alpha+\beta}=\sqrt{\la}T\wt{x}_1 \ \ \mbox{ or } \ \ \alpha+\frac\beta 2 = \la T^2\wt{x}_1^2 (\alpha+\beta)^2.
\end{equation}
Thus we need to minimise $\sqrt{\alpha + \beta/2}$ on the set $\alpha, \beta >0$ knowing that the condition \reff{10.22} holds true. Let $c=\lambda T^2\wt{x}_1^2$. For positive $\beta, \alpha$ the condition \reff{10.22} draws the curve $\mathcal{C}$
$$
\mathcal{C} = \Bigl\{ (\alpha, \beta):\  \alpha = \frac{1}{2c} - \beta + \sqrt{ 1-2c\beta},\ \ \beta \in \Bigr[0, \frac{1}{2c} \Bigr] \Bigr\}.
$$
The minimum of $\sqrt{\alpha+\frac\beta 2}$ is attained at the same point as the minimum of $\alpha+\frac\beta 2$. Thus, on the set $\mathcal{C}$, we have 
$$
\alpha+\frac\beta 2 = \frac{1}{2c} - \frac\beta 2 + \sqrt{ 1-2c\beta}.  
$$
It means that the minimum is attained at the point $\beta = \hat\beta := \bigl(2\lambda T^2\wt{x}_1^2 \bigr)^{-1}$ which corresponds to $\alpha=0$. 
Therefore, the minimal value of the rate function equals to
\[
J(0,\hat\beta)=2+ \frac1{2\wt{x}_1}.
\]
\qed

\section*{Acknowledgement}
The researches of E. Pechersky, S. Pirogov and A. Vladimirov were carried out
at the Institute for Information Transmission Problems (IITP), Russian
Academy of Science and funded by the Russian Science Foundation  (grant
14-50-00150). 
The main contribution by S. Pirogov was done in the section 5. The main contribution by A. Vladimirov was done in the section 6. The main contribution by E. Pechersky was done in the section 7.

G.M. Sch\"utz thanks University of S\~ao Paulo (USP) for kind hospitality.

A. Yambartsev   thanks  Conselho Nacional de Desenvolvimento
Cient\'ifico e Tecnol\^ogico (CNPq) grant 301050/2016-3 and Funda\c{c}\~ao de Amparo \`a Pesquisa do Estado de S\~ao Paulo (FAPESP) grant 2017/10555-0.



\end{document}